\documentclass[a4paper,10pt]{amsart}

\usepackage{amsfonts,amsthm,amssymb,amsmath,amscd}

\newcommand{\ZZ}{\mathbb{Z}}
\newcommand{\QQ}{\mathbb{Q}}
\newcommand{\RR}{\mathbb{R}}
\newcommand{\CC}{\mathbb{C}}
\newcommand{\PP}{\mathbb{P}}

\newcommand{\la}{\lambda}

\newcommand{\DD}{\mathbb{D}}
\newcommand{\HH}{\mathbb{H}}
\newcommand{\Si}{\mathfrak{S}}

\newcommand{\GL}{\mathrm{GL}}
\newcommand{\SL}{\mathrm{SL}}

\newcommand{\G}{\Gamma}

\newcommand{\I}{\mathrm{I}}
\newcommand{\IM}{\mathrm{Im}}

\newcommand{\Ons}{\mathrm{O}^+_{ns}}

\begin{document}
\title{Hessian K3 surfaces of non Sylvester type}
\date{}
\author{Kenji Koike, Yamanashi University}  
\maketitle
\begin{abstract}
We construct the moduli space of cubic surfaces which do not admit a Sylvester form as 
an arithmetic quotient, and determine the graded ring of modular forms of even weights.
\end{abstract} 
\section{Hessian K3 surfaces of non Sylvester type}
\subsection{}
It is classicaly known that the ring of $\SL_4(\CC)$-invariants of quaternary cubic forms is  
\[
  \CC[I_8, I_{16}, I_{24}, I_{32}, I_{40}, I_{100}] \qquad (\deg I_n = n)
\] 
where $I_8, \cdots, I_{40}$ are algebraically independent and $I_{100}^2 \in \CC[I_8, I_{16}, I_{24}, I_{32}, I_{40}]$ 
(\cite{H}, \cite{Sa}). Hence the moduli space of cubic surfaces $\mathcal{M}_I$ is isomorphic to the weighted projective space 
\[
 \mathrm{Proj} \medspace \CC[I_8, I_{16}, I_{24}, I_{32}, I_{40}] = \PP(1,2,3,4,5)_I.
\] 
A general cubic surface is written as a complete intersection
\[
 S_{\la} \ : \ X_0 + \cdots + X_4 = 0, \quad 
\la_0 X_0^3 + \cdots + \la_4 X_4^3 = 0 
\]
in $\PP^4$ with $\la_0, \cdots, \la_4 \ne 0$, which is called the Sylvester form. Let $\sigma_i$ be 
the $i$-th elementary symmetric polynomial in $\la_0, \cdots, \la_4$. They give invariants of $S_{\la}$, and 
we have
\begin{align*}
 I_8 = \sigma_4^2 - 4 \sigma_3 \sigma_5, \quad I_{16} = \sigma_5^3 \sigma_1, \quad
I_{24} = \sigma_5^4 \sigma_4, \quad I_{32} =\sigma_5^6 \sigma_2, \quad I_{40} = \sigma_5^8.
\end{align*}
This correspondence gives a birational map
\[
 \PP(1,2,3,4,5)_{\la} \longrightarrow \PP(1,2,3,4,5)_I
\]
with the base locus $\sigma_5 = \sigma_4 = 0$.
The Hessian of $S_{\la}$ is given by
\[
 H_{\la} \ : \ \ X_0 + \cdots + X_4 = 0, \quad 
\frac{1}{\la_0 X_0} + \cdots + \frac{1}{\la_4 X_4} = 0.
\]
The Picard lattice of the desingularization of a general $H_{\la}$ is 
$\mathrm{U} \oplus \mathrm{U}(2) \oplus \mathrm{A}_2(2)$ (see \cite{DK}).
\subsection{Dardanelli - van Geemen's stratification}
The following facts on $\mathcal{M}$ were proved in \cite{DvG}.
\\
(I) \ The subvariety of $\mathcal{M}$ parametrizing cubic surfaces which do not admit a Sylvester form is defined by $I_{40} = 0$. 
In general, such surfaces are given by
\[
 S_{ns1}(a) : X_1^3 + X_2^3 + X_3^3 - X_0^2 (a_0 X_0 + 3 a_1 X_1 + 3 a_2 X_2 + 3 a_3 X_3) = 0.
\]
If we denote the $i$-th elementary symmetric polynomial in $a_1^3, \ a_2^3,\ a_3^3$ by $\rho_i$, then we have
\[
 [S_{ns1}(a)] = [-4 \rho_1 + a_0^2:\rho_2:2\rho_3:\rho_1 \rho_3:0] \in \PP(1,2,3,4,5)_I.
\]
The Hessian surface of $S_{ns1}(a)$ is given by
\[
 H_{ns1}(a) : X_0 X_1 X_2 X_3(a_1 \frac{X_1}{X_0} + a_2 \frac{X_2}{X_0} + a_3 \frac{X_3}{X_0} + a_0
+ a_1^2 \frac{X_0}{X_1} + a_2^2 \frac{X_0}{X_2} + a_3^2 \frac{X_0}{X_3}) = 0,
\]
and the transcendental lattice of the desingularization of a general $H_{ns1}(a)$ is 
$T_{ns1} = \mathrm{U} \oplus \mathrm{U}(2) \oplus \left< -4 \right>$. 
In affine coordinates $[X_0:X_1:X_2:X_3] = [1:x/a_1:y/a_2:z/a_3]$, the equation of $H_{ns1}(a)$ is 
\[
 xyz(x + y + z + a_0 + a_1^3 \frac{1}{x} + a_2^3 \frac{1}{y} + a_3^3 \frac{1}{z}) =0.
\]
(II) \ The subvariety of $\mathcal{M}$ parametrizing cubic surfaces
\[
 S_{ns2}(b) : X_1^3 + X_2^3 + 2 b_0 X_3^3 - 3X_3(b_1 X_1 X_3 + X_2 X_3 + X_0^2 ) = 0
\]
is defined by $I_{24} = I_{40} = 0$, and we have
\[
 [S_{ns2}(b)] = [-8b_0:1+b_1^3:0:b_1^3:0] \in \PP(1,2,3,4,5)_I.
\]
The Hessian surface of $S_{ns2}(b)$ is given by
\[
 H_{ns2}(b) : X_1 X_2 X_3 (-2b_0 X_3 + b_1 X_1 + X_2) + X_3^3 (X_1 + b_1^2 X_2) - X_0^2 X_1 X_2 =0
\]
and the transcendental lattice of the desingularization of a general $H_{ns2}(b)$ is 
$T_{ns2} = \mathrm{U} \oplus \mathrm{U}(2)$.
\\
(III) \ The subvariety of $\mathcal{M}$ parametrizing ``cyclic cubic surfaces''
\[
 S_{cyc}(a) : a_4 X_4^3 - a_3(X_0 + X_1 + X_2)^3 + a_0 X_0^3 + a_1 X_1^3 + a_2 X_2^3 = 0
\]
is defined by $I_{24} = I_{32} = I_{40} = 0$, and we have
\[
 [S_{cyc}(a)] \in [\mu_3^2 - 4 \mu_2 \mu_4 : \mu_4^3:0:0:0] \in \PP(1,2,3,4,5)_I
\]
where $\mu_i$ is the $i$-th symmetric polynomial of $a_0, \cdots, a_3$. The Hessian of $S_{cyc}(a)$ is reducible.
\\
(IV) \ The strictly semi-stable surface $t^3 = xyz$ corresponds to the point $[8:1:0:0:0]$, and the Fermat cubic surface 
corresponds to the point $[1:0:0:0:0]$.
\subsection{Batyrev's mirror construction} 
Hessian surfaces $\{H_{ns1}\}$ are obtained also as toric hypersurfaces.  
Let $\Delta$ be the octahedron in $\RR^3$ with vertices 
\[
 (\pm1, 0, 0), \quad (0, \pm 1, 0), \quad (0, 0, \pm1).
\]
It is a simplicial reflexive polytope, and its dual polytope $\Delta^*$ is the cube with vertices $(\pm1, \pm1, \pm1)$. 
Considering faces of $\Delta$ as simplicial cones, we obtain a toric variety $X(\Delta) = \PP^1 \times \PP^1 \times \PP^1$. 
The linear system of anti-canonical classes of $X(\Delta)$ (that is, K3 surfaces of degree $(2,2,2)$ in 
$(\PP^1)^3$) is given by 
\[
\mathcal{F}(\Delta^*) = \{\sum a_{ijk} x^i y^j z^k \ = 0 | \ (i,j,k) \in \Delta^* \cap \ZZ^3 \} \quad
((x,y,z) \in (\CC^{\times})^3 \subset (\PP^1)^3).
\]
Similarly, we have the dual family of K3 surfaces 
\[
\mathcal{F}(\Delta) = \{c_1x + c_2y + c_3z +c_4 + c_5 \frac{1}{x} + c_6 \frac{1}{y} + c_7\frac{1}{z} = 0\}
\]
as hypersurfaces of $X(\Delta^*)$. It is obvious that this family is birationally equivalent to the family $\{H_{ns1}\}$. 
Note that the Picard lattice of a general menber of $\mathcal{F}(\Delta^*)$ is 
\[
 P = \begin{bmatrix} 0 & 2 & 2 \\ 2 & 0 & 2 \\ 2 & 2 & 0 \end{bmatrix} 
\cong \mathrm{U}(2) \oplus \left< -4 \right>,
\]
and we have $T_{ns1} = \mathrm{U} \oplus P$. Hence $\mathcal{F}(\Delta)$ is the mirror partner 
of $\mathcal{F}(\Delta^*)$ (see \cite{B}, \cite{C1}, \cite{D1} and \cite{GN}). Note also that $\mathcal{F}(\Delta)$ 
is a subfamily of $\mathcal{F}(\Delta^*)$. In the following, we regard $H_{ns1}$ as hypersurfaces in 
$(\PP^1)^3$, and  we replace coefficients $a_0, a_1^3, a_2^3, a_3^3$ of $H_{ns1}$ by 
$1, u_1, u_2, u_3$: 
\[
 H(u) : f_u = xyz(x + y + z + 1) + (u_1 y z + u_2 z x + u_3 x y) = 0 \qquad (x,y,z) \in (\PP^1)^3.
\]
\subsection{Remark}
From the $1$-parameter family 
\[
 H_{PS}(u) : xyz(x+y+z+1) + u(xy + yz + zx) = 0,
\]
by the base change $u = (t + t^{-1})^{-2}$, we obtain the family
\[
 x + \frac{1}{x} + y + \frac{1}{y} + z + \frac{1}{z} + t + \frac{1}{t} =0
\]
studied by Peters and Stienstrain in \cite{PS}. They studied the Picard-Fuchs equation and modular forms. 
The transcendental lattice of a general member is $\mathrm{U} \oplus \left< 12 \right>$. 
This K3-fibration is considered as a (singular) Calabi-Yau hypersurface in $(\PP^1)^4$ (see \cite{V}).
\subsection{Singularities}
Let us assume $u_1 u_2 u_3 \ne 0$. Then $H(u) \cap (\CC^{\times})^3$ is smooth 
if and only if
\[
  \Delta_{sing}(u) = \prod (1 \pm 2\sqrt{u_1} \pm 2\sqrt{u_2} \pm 2\sqrt{u_3}) \ne 0.
\]
Therefore we define the parameter space 
\[
 \mathcal{U} = \{ u = (u_1, u_2, u_3) \ | \ u_1 u_2 u_3 \Delta_{sing}(u) \ne 0 \}.
\]
For any $u \in \mathcal{U}$, we see that $H(u) \cap ((\PP^1)^3 - (\CC^{\times})^3)$ is decomposed into twelve lines
\begin{align*}
L_{x00} = \PP^1 \times \{0\} \times \{0\}, \quad L_{x0\infty} = \PP^1 \times \{0\} \times \{\infty\},\ \cdots,
\ L_{\infty \infty z} = \{\infty\} \times \{\infty\} \times \PP^1.
\end{align*}
They intersect at eight points
\begin{align*}
 (0,0,0), \quad (0,0,\infty), \quad (0,\infty,0), \quad (0,\infty,\infty),\\
(\infty,0,0), \quad (\infty,0,\infty), \quad (\infty,\infty,0), \quad (\infty,\infty,\infty),
\end{align*}
that are singular points of $H(u)$, and all of them are $A_1$-singularities. Blowing up eight singular points of $H(u)$, 
we obtain a K3 surface $\widetilde{H}(u)$. Let $N_u \subset \mathrm{H}^2(\widetilde{H}(u), \ZZ)$ be a sublattice generated by twelve lines 
$L_{x00}, \cdots, L_{\infty \infty z}$ and eight exceptional curves $E_{000}, \cdots, E_{\infty \infty \infty}$ that are blown down to 
$(0,0,0), \cdots, (\infty \infty \infty)$.
\subsection{Proposition} 
\label{prop-inv}
(1) \ For a general $u \in \mathcal{U}$, the lattice $N_u$ is the Picard lattice $\mathrm{Pic}(\widetilde{H}(u))$.
\\
(2) \ We have three involutions
\[ 
\epsilon_x : (x,y,z) \mapsto (\frac{u_1}{x}, y, z), \quad
\epsilon_y : (x,y,z) \mapsto (x, \frac{u_2}{y}, z), \quad
\epsilon_z : (x,y,z) \mapsto (x, y, \frac{u_3}{z}).
\]
on $\widetilde{H}(u)$, and the product $\epsilon = \epsilon_x \epsilon_y \epsilon_z$ is an Enriques involution.
\\
(3) Let $N_u^* \subset N \otimes \QQ$ be the dual lattice of $N_u$, and 
$q_N : N_u^* / N_u \rightarrow \QQ / 2 \ZZ$ be the discriminant form (\cite{N}). Then we have 
$\epsilon = \epsilon_x = \epsilon_y = \epsilon_z$ as elements of the finite orthogonal group $\mathrm{O}(q_N)$. 
Moreover, we have $\mathrm{O}(q_N) = \mathrm{S}_3 \times \left< \epsilon \right>$, where $\mathrm{S}_3$ 
is realized as symmetry of $(x,y,z)$.
\\ \\
{\bf Proof.} 
(1) \ The self intersection numbers of $L_{***}$ and $E_{***}$ are $-2$, and we have $E_{abc} \cdot L_{stu} = 1$ if
two of three equalities $a=s$ , $b=t$ or $c=u$ are hold. Other intersection numbers are zero. Using a computer, we can show that the rank of 
the intersection matrix of them is $17$. In fact, we have equalities 
\begin{align*}
 E_{000} = E_{00 \infty} + E_{0\infty 0} + 3E_{0\infty \infty} -3E_{\infty 00} -  E_{\infty 0 \infty} - E_{\infty \infty 0} 
+ E_{\infty \infty \infty} \\
-2 L_{x00} + 2L_{x \infty \infty}+ 2L_{0y \infty} -2 L_{\infty y0} +2 L_{0\infty z} -2 L_{\infty 0z}, 
\\
L_{\infty y \infty} = 2E_{0\infty 0} + 2E_{0\infty \infty} -2E_{\infty 00} - 2E_{\infty 0 \infty} - L_{x00} - L_{x0 \infty} 
+ L_{x \infty 0} \\ 
+ L_{x \infty \infty} + L_{0y0} + L_{0y \infty} - L_{\infty y 0} + 2L_{0 \infty z}  -2 L_{\infty 0 z},
\\
L_{\infty \infty z} = 2E_{00 \infty} + 2E_{0\infty \infty} -2E_{\infty 00} - 2E_{\infty \infty 0} - L_{x00} + L_{x0 \infty} 
- L_{x \infty 0} \\
+ L_{x \infty \infty} + 2L_{0y \infty} -2 L_{\infty y0} + L_{0 \infty z}  - L_{\infty 0 z}
\end{align*}
as elements of $N_u$. Therefore $E_{000}$, $L_{\infty y \infty}$ and $L_{\infty \infty z}$ are redundant. 
Since the determinat of the intersection matrix of other $17$ curves is $16$, we see that they span the orthogonal complement of
$T_{ns1} = \mathrm{U} \oplus \mathrm{U}(2) \oplus \left< -4 \right>$.
\\
(2) \ As an involution of $(\PP^1)^3$, fixed points of $\epsilon$ are 
$(\pm \sqrt{u_1}, \pm \sqrt{u_2}, \pm \sqrt{u_3})$. If $u \in \mathcal{U}$, then such points are not on $H(u)$. 
\\
(3) \ We have $N_u^* / N_u \cong T_{ns1}^* / T_{ns1} \cong (\ZZ / 2 \ZZ)^2 \times (\ZZ / 4\ZZ)$, and it is generated by
\begin{align*}
\ell_1 =& \frac{1}{2}(L_{0y0} + L_{0y \infty} + L_{00z} + L_{0 \infty z}), \\
\ell_2 =& \frac{1}{2}(L_{x00} + L_{x0 \infty} + L_{00z} + L_{\infty 0 z}), \\
m =& \frac{1}{4}(2 E_{0\infty \infty} + 2E_{\infty0\infty} + 2E_{\infty \infty0} + 2E_{\infty \infty \infty} + 2L_{x00} \\
&+ 3L_{x0\infty} + 3L_{x\infty0} + 2L_{0y0} + L_{0y\infty} + L_{\infty y0} + 3L_{0\infty z} + L_{\infty 0z}).
\end{align*}
By machine computation, we see that 
\[
 \epsilon_x(\ell_i) = \ell_i \ (i=1,2), \qquad \epsilon_x(m) = -m, 
\]
and the same for $\epsilon_y$ and $\epsilon_z$. The $2$-torsion subgroup of $N_u^* / N_u$ is generated by
$\ell_1, \ell_2$ and $\ell_3 = 2m + \ell_1 + \ell_2$, and these are all of elements $x \in N_u^* / N_u$ of order $2$ such that $q_N(x) = 0$. 
We have a split exact sequence
\[
 1 \longrightarrow \left< \epsilon \right> \longrightarrow \mathrm{O}(q_N) \longrightarrow 
\{\text{permutations of $\ell_1, \ell_2, \ell_3$}\} \longrightarrow 1  
\] 
since permutations of $(x,y,z)$ give permutations of $\ell_i$'s.
\hfill $\Box$
\section{The period mapping and modular groups} 
\subsection{The period mapping}
The period domain of the family $\{\widetilde{H}(u) \ | \ u \in \mathcal{U}\}$ is the bounded symmetric domain 
\begin{align*}
\DD_{ns} = \{z \in \PP^4 \ |\ {}^tzQz = 0, \ {}^tzQ\bar{z} > 0 \}, \quad 
Q = \begin{bmatrix} 0 & 1 \\ 1 & 0 \end{bmatrix} \oplus \begin{bmatrix} 0 & 2 \\ 2 & 0 \end{bmatrix} \oplus [-4].
\end{align*}
of type IV defined by the lattice $T_{ns1}$. More explicitly, we have
\begin{align*}
[1: z_2: \cdots : z_5] \in \DD_{ns} \Leftrightarrow \begin{cases} 
z_2 = -2 (z_3 z_4 - z_5^2) \\
y_3 y_4 - y_5^2 >0 \ \ (y_i = \mathrm{Im} z_i) \end{cases}
\end{align*}
and $\DD_{ns} = \DD_{ns}^+ \coprod \DD_{ns}^-$ where $\DD_{ns}^{\pm} = \{z \in \DD_{ns} \ : \ \pm y_3 > 0 \}$. 
Let us define the orthogonal group
\begin{align*}
\Ons = \{ g \in \GL_5(\ZZ) \ | \ {}^tgQg = Q, \ g(\DD_{ns}^+) = \DD_{ns}^+ \}
\end{align*}
on the lattice $T_{ns1}$, which acts on $\DD_{ns}^+$. We define also the discriminant form
\begin{align*}
q_{ns1} : T_{ns1}^* / T_{ns1} \longrightarrow \QQ / 2 \ZZ
\end{align*} 
and the orthogonal group $\mathrm{O}(q_{ns1})$. Let $\Ons(2)_{\epsilon}$ be the kernel of the natural homomorphism 
\[
 \Ons \longrightarrow \mathrm{O}(q_{ns1}) \cong \mathrm{O}(q_{N}) \cong \mathrm{S}_3 \times \left< \epsilon \right>,
\]
and $\Ons(2)$ be the kernel of the composition map
\[
 \Ons(2)_{\epsilon} \longrightarrow \mathrm{S}_3 \times \left< \epsilon \right> \longrightarrow \mathrm{S}_3. 
\]
We have $-1 \in \Ons(2)$ and $-1 \notin \Ons(2)_{\epsilon}$. Since $[\Ons(2)_{\epsilon} : \Ons(2)] = 2$, we see that
\[
 \DD_{ns}^+ / \Ons(2)_{\epsilon} = \DD_{ns}^+ / \Ons(2).
\]
Let $S_u \subset \mathrm{H}_2(\widetilde{H}(u), \ZZ)$ be the sublattice generated by $L_{***}$'s and $E_{***}$'s, that is, 
the Poincare dual of $N_u \subset \mathrm{H}^2(\widetilde{H}(u), \ZZ)$. Taking suitable $2$-cycles
$\gamma_1(u), \cdots, \gamma_5(u) \in (S_u)^{\perp} \cong T_{ns1}$ that are uniquely determined  
up to $\Ons$-action, we can define the period mapping
\[
 Per : \mathcal{U} \longrightarrow  \DD_{ns}^+, \quad 
u = (u_1,u_2,u_3) \mapsto [\int_{\gamma_1(u)} \omega_u : \cdots : \int_{\gamma_5(u)} \omega_u]
\]
where $\omega_u \in \mathrm{H}^{2,0}(\widetilde{H}(u))$.
\subsection{Proposition}
The multi-valued  map $Per$ induces an injective $\mathrm{S}_3$-equivariant map 
$\mathcal{U} \longrightarrow  \DD_{ns}^+ / \Ons(2)_{\epsilon}$ and the map
$\mathcal{U} / \mathrm{S}_3 \longrightarrow  \DD_{ns}^+ / \Ons$ for $\mathrm{S}_3$-quotients.
\\ \\
{\bf Proof.} 
Note that \\  
(1) the monodromy action of $\pi_1(\mathcal{U}, u)$ on $S_u \subset \mathrm{H}_2(\widetilde{H}(u), \ZZ)$ is trivial, 
\\
(2) we can lift $g \in \Ons$ to $\tilde{g} \in \mathrm{O}(\mathrm{H}_2(\widetilde{H}(u), \ZZ))$ such that 
$\tilde{g}|_{S_u} = \mathrm{id}$ iff $g \in \Ons(2)_{\epsilon}$.
\\
From these facts together with Proposition \ref{prop-inv}, we see that the map is injective
as the period map of $N_u$-polarized K3-surfaces (see \cite{D1}).
\hfill $\Box$
\subsection{Proposition} \label{FC-prop}
The period map $Per$ is given by the developping map of the Lauricella's hypergeometric differential equation for 
$F_C(1, \frac{1}{2} ; 1,1,1 ; -2u_1, -2u_2, -2u_3)$ (see \cite{Y}).
\\ \\
{\bf Proof.}
Indeed, we obtain a period of $H(u)$ as follows.
\begin{align*}
I(u_1,u_2,u_3) &= \iiint_{|x|=|y|=|z|=\varepsilon} 
\frac{dx \wedge dy \wedge dz}{f_u} \\
&= \iiint_{|x|=|y|=|z|=\varepsilon} 
\frac{1}{xyz(x + y + z + 1)} \frac{dx dy dz}{1 + \frac{u_1 xy + u_2 yz + u_3 zx}{xyz(x + y + z + 1)}} \\
&= \iiint_{|x|=|y|=|z|=\varepsilon} 
\sum_{n=0}^{\infty} \frac{(-u_1 xy - u_2 yz -u_3 zx)^n}{(xyz(x + y + z + 1))^{n+1}} dx dy dz 
\quad (|u_i| \ll \varepsilon) \\
&= \iiint_{|x|=|y|=|z|=\varepsilon}
 \sum_{p,q,r=0}^{\infty} \frac{(p+q+r)!}{p!q!r!} 
\frac{x^{p+r} y^{p+q} z^{q+r} dx dy dz }{(xyz(x + y + z + 1))^{p+q+r+1}} (-u_1)^p (-u_2)^q (-u_3)^r \\
&= \sum_{p,q,r=0}^{\infty} \frac{(p+q+r)!}{p!q!r!} N(p,q,r)(-u_1)^p (-u_2)^q (-u_3)^r 
\end{align*}
where
\begin{align*}
 N(p,q,r) &= \iiint_{|x|=|y|=|z|=\varepsilon} \frac{dx dy dz }{x^{q+1} y^{r+1} z^{p+1} (x + y + z + 1)^{p+q+r+1}} \\
&= (2 \pi i)^3 \frac{(2p+2q+2r)!}{(p+q+r)!p!q!r!}.
\end{align*}
Therefore, we obtain 
\begin{align*}
 I(u_1,u_2,u_3) &= (2 \pi i)^3 \sum_{p,q,r=0}^{\infty} \frac{(2p+2q+2r)!}{(p!q!r!)^2} (-u_1)^p (-u_2)^q (-u_3)^r \\
&= \sum_{p,q,r=0}^{\infty} \frac{(1)_{p+q+r} (\frac{1}{2})_{p+q+r}}{(1)_p (1)_p (1)_q (1)_q (1)_r (1)_r} 
(-2u_1)^p (-2u_2)^q (-2u_3)^r \\
&= F_C(1, \frac{1}{2} ; 1,1,1 ; -2u_1, -2u_2, -2u_3).
\end{align*}
\hfill $\Box$
\subsection{Modular groups}
The domain $\DD_{ns}^+$ is isomorphic to the Siegel upper half space $\Si_2$ of degree $2$ by the map
\[
 \Psi :  \DD_{ns}^+ \longrightarrow \Si_2 = \{ \tau \in \GL_2(\CC) \ | \ \IM \tau > 0  \}, \qquad 
[1:z_2:\cdots:z_5] \mapsto \begin{bmatrix} z_3 & z_5 \\ z_5 & z_4 \end{bmatrix}.
\]
The symplectic group 
\[
 \mathrm{Sp}_{2g}(\RR) = \{ g \in \GL_{2g}(\RR) \ | \ {}^tg J g = J \}, \quad J = \begin{bmatrix} 0 & -\I_g \\ \I_g & 0 \end{bmatrix}
\]
acts on $\Si_g$ by $\begin{bmatrix} A & B \\ C & D \end{bmatrix} \cdot \tau = (A \tau +B)(C \tau +D)^{-1}$.
Let $\G_g$ be the Siegel modular group $\mathrm{Sp}_{2g}(\RR) \cap \GL_{2g}(\ZZ)$. We consider the congruence subgroup
\begin{align*}
\G_0(2)_g = \{ \begin{bmatrix} A & B \\ C & D \end{bmatrix} \in \G_g \ | \ C \equiv 0 \mod 2 \}, \quad
\end{align*}
and the extension $\G_0^*(2)_2$ of $\G_0(2)_2$ by a normalizer $W = \displaystyle
\frac{1}{\sqrt{2}} \begin{bmatrix} 0 & -\I_2 \\ 2 \I_2 & 0\end{bmatrix}$.
\subsection{Proposition}
 Then we have an isomorphism 
$\Ons / \{\pm1 \} \cong \G_0^*(2)_2 / \{\pm 1\}$ as aotumorphisms of 
$\DD_{ns}^+ \cong \Si_2$.
\\ \\
{\bf Proof.} 
This is an easy consequence of Theorem 3.1 in \cite{Ko2}, and we omit the proof. 
We give just explicit correspondences of generators: 
\\
(1) The map $g : \GL_2(\ZZ) \rightarrow \Ons$,
\[ 
\begin{bmatrix} a_1 & a_2 \\ a_3 & a_4 \end{bmatrix} \mapsto 
\I_2 \oplus \begin{bmatrix} 
a_1^2 & a_2^2 & 2 a_1 a_2 \\
a_3^2 & a_4^2 & 2 a_3 a_4 \\ 
a_1 a_3 & a_2 a_4 & a_1 a_4 + a_2 a_3 
\end{bmatrix} 
\]
is a homomorphism such that $\mathrm{Ker} \ g = \{ \pm1 \}$ and 
$\Psi(g(A) \cdot z) = \begin{bmatrix} A & 0 \\ 0 & {}^tA^{-1} \end{bmatrix} \cdot \Psi(z)$. 
\\
(2) Let $\mathcal{B}_2$ be the additive group of integral symmetric matrices of degree $2$. Then 
the map $h : \mathcal{B}_2 \rightarrow \Ons$,
\begin{align*}
\begin{bmatrix} m_1 & m_2 \\ m_2 & m_3 \end{bmatrix}
\mapsto \begin{bmatrix} 1 & 0 & 0 & 0 & 0 \\
-2m_1 m_2 + 2 m_3^2 & 1 & -2m_2 & -2m_1 & 4m_3 \\ 
m_1 & 0 & 1 & 0 & 0\\
m_2 & 0 & 0 & 1 & 0\\
m_3 & 0 & 0 & 0 & 1\\
\end{bmatrix}
\end{align*}
is a homomorphism such that $\Psi(h(B) \cdot z) = \begin{bmatrix} \I_2 & B \\ 0 & \I_2 \end{bmatrix} \cdot \Psi(z)$.
\\
(3) For $w = \begin{bmatrix} 0 & 1 \\ 1 & 0 \end{bmatrix} \oplus \begin{bmatrix} 0 & 1 \\ 1 & 0 \end{bmatrix} 
\oplus [-1] \in \Ons$, we have
$\displaystyle \Psi(w \cdot z) = -\frac{1}{2} \Psi(z)^{-1} = W \cdot \Psi(z)$. 
\hfill $\Box$
\subsection{Proposition} \label{prop-table}
(1) If $x, y \in \ZZ^5$ satisfy ${}^tx Qx = {}^tyQy = 0$ and ${}^txQy = 1$, then there exists a transformation 
$\gamma \in \Ons$ such that $\gamma \cdot x = e_1$ and $\gamma \cdot y = e_2$, where $e_i$ is the $i$-th unit 
vector.
\\
(2) For any primitive sublattice $M \cong \mathrm{U} \oplus \left< 12 \right>$ of $T_{ns1}$, 
there exists $\gamma \in \Ons$ such that $\gamma(M)$ is either
\[
 M_1 = \ZZ e_1 \oplus \ZZ e_2 \oplus \ZZ (e_3 + 3 e_4) \quad \text{or} \quad
 M_2 = \ZZ e_1 \oplus \ZZ e_2 \oplus \ZZ (2 e_3 + 2 e_4 + e_5).
\]
For any primitive sublattice $M' \cong \mathrm{U} \oplus \mathrm{U}(2)$ of $T_{ns1}$, 
there exists $\gamma' \in \Ons$ such that 
\[
 \gamma'(M') = \ZZ e_1 \oplus \ZZ e_2 \oplus \ZZ e_3 \oplus \ZZ e_4.
\]
(3) We have the following table for periods of special subfamilies:
\\ \\
\begin{tabular}{c|c|c|c|c}
& $\PP(1:2:3:4:5)_I$ & lattice & $\DD_{ns}^+ / \Ons$ & $\Si_2 / \G_0^*(2)_2$  
\\
\hline
$H_{ns2}(b)$ & $[-8b_0:1+b_1^3:0:b_1^3:0]$ & $\mathrm{U} \oplus \mathrm{U}(2)$ & $[1:z_2:z_3:z_4:0]$ & 
$\begin{bmatrix} z_3 & 0 \\ 0 & z_4 \end{bmatrix}$
\\
\hline
$H_{PS}(u)$ & $[-12u +1:3u^2:2u^3:3u^4:0]$ & $\mathrm{U} \oplus \left< 12 \right>$ & $[1:z_2:2z_5:2z_5:z_5]$ & 
$\begin{bmatrix} 2z_5 & z_5 \\ z_5 & 2z_5 \end{bmatrix}$
\end{tabular}
\\ \\
{\bf Proof.}
(1) This is shown by the same argument with Proposition 3.2 in \cite{Ko2}.
\\
(2) By (1), there exists $\gamma \in \Ons$ and $x,y,z \in \ZZ$ such that
\[
 \gamma(M) = \ZZ e_1 \oplus \ZZ e_2 \oplus \ZZ(xe_3 + y e_4 + z e_5), \quad xy -z^2 = 3.
\]
Now the assertion for $M$ follows from the facts:
\\
(i) the integer solutions of the system of eqations
\[
 xy = z^2 + 3, \qquad |x|>|z|, \qquad |y|>|z|
\]
are $(2,2,\pm 1)$ or $(-2,-2,\pm 1)$,
\\
(ii) if $|x|<|z|$ or $|y|<|z|$, then multiplying
\[
 \I_2 \oplus \begin{bmatrix} 1 & 1 & \pm 2 \\ 0 & 1 & 0 \\ 0 & \pm 1 & 1 \end{bmatrix}, \qquad
\I_2 \oplus \begin{bmatrix} 1 & 0 & 0 \\ 1 & 1 & \pm 2 \\ \pm 1 & 0 & 1 \end{bmatrix} \in \Ons,
\]
we can decrease the value of $|z|$. 
\\
The assertion for $M'$ is easily shown by the same way.
\\
(3) By (2), periods of $H_{ns2}(b)$ belong to the divisor $\{z_5 = 0\}$ in $\DD_{ns}^+$. 
Because surfaces $H_{PS}(u)$ don't belong to the family $\{H_{ns2}(b)\}$, their periods don't belong to 
$\PP(M_1 \otimes \CC) \subset \{z_5 = 0\}$. Hence periods of $H_{PS}(u)$ belong to $\PP(M_2 \otimes \CC)$. 
\hfill $\Box$
\section{Graded ring of theta constants}
\subsection{}
Let $\G'$ be a subgroup of $\mathrm{Sp}_4(\RR)$. A holomorphic function $f(\tau)$ on $\Si_2$ is a modular form of weight $k$ 
with respect to $\G'$ if it holds
\[
 f((A \tau + B)(C \tau + D)^{-1}) = \det (C \tau + D)^k f(\tau) 
\]
for any $\begin{bmatrix} A & B \\ C & D \end{bmatrix} \in \G'$. Let $M_k(\G')$ be the vector space of 
such functions, and $A(\G')_{even}$ be the graded ring $\oplus_{k=0}^{\infty} M_{2k}(\G')$.
The generaters of the graded ring $A(\G_0(2)_2)_{even}$ are given by theta constants
\[
 \theta_{a,b}(\tau) = \sum_{n \in \ZZ^2} \exp[\pi i {}^t(n+a)\tau(n+a) + 2 \pi i {}^t(n+a)b], \quad \tau \in \Si.
\]
For simplicity, we denote $\theta_{a,b}$ by $\theta_{xyzw}$ if $a = {}^t(x/2, y/2)$ and $b = {}^t(z/2, w/2)$. 
\subsection{Theorem{\rm (Ibukiyama, \cite{Ib})}} 
Let us define modular forms
\begin{align*}
\vartheta = (\theta_{0000}^4 + \theta_{0001}^4 + \theta_{0010}^4 + \theta_{0011}^4)/4, \quad
\phi_1 = (\theta_{0000} \theta_{0001} \theta_{0010} \theta_{0011})^2, \quad \phi_2 &= (\theta_{0100}^4 - \theta_{0110}^4)^2 / 16384, \\
\chi = (\theta_{0100} \theta_{0110} \theta_{1000} \theta_{1001} \theta_{1100} \theta_{1111})^2 / 4096
\end{align*}
of weight $2, 4, 4$ and $6$. 
Then the graded ring $A(\G_0(2)_2)_{even}$ is a free algebra $\CC[\vartheta, \phi_1, \phi_2, \chi]$, and
\[
 \mathrm{Proj}\ A(\G_0(2)_2)_{even} \cong \PP(2,4,4,6).
\]
\subsection{Lemma} \label{chi-vanish}
The zero divisor of the function $\chi (\tau)$ is $\G_0(2)_2$-obit of 
\[
 \HH \times \HH = \{ \tau = \begin{bmatrix} \tau_1 & \tau_2 \\ \tau_2 & \tau_3 \end{bmatrix} \in \Si\ | \ \tau_2 = 0 \}
\]
with multiplicity $1$, and $\chi(\tau)$ is the unique non-trivial function in $M_6(\G_0(2)_2)$ vanishing there. 
\\ \\
{\bf Proof.}
The first assertion is proved by exactly the same way as in \cite{Kl}, p.116 - p.118.
By the equality of theta constants of one variable $\theta_{00}^4 = \theta_{01}^4 + \theta_{10}^4$, we see that
\begin{align*}
\vartheta(\tau) &= (\theta_{00}^4(\tau_1) + \theta_{01}^4(\tau_1))(\theta_{00}^4(\tau_3) + \theta_{01}^4(\tau_3))/4, \\ 
\phi_1 &= \theta_{00}^4(\tau_1) \theta_{01}^4(\tau_1) \theta_{00}^4(\tau_3) \theta_{01}^4(\tau_3), \\
\phi_2 &= (\theta_{00}^4(\tau_1) - \theta_{01}^4(\tau_1))^2 (\theta_{00}^4(\tau_3) - \theta_{01}^4(\tau_3))^2 / 16384
\end{align*}
for $\tau \in \HH \times \HH$. Therefore $\vartheta^3, \ \vartheta \phi_1, \ \vartheta \phi_2$ are linearly independent on $\HH \times \HH$.
\hfill $\Box$
\subsection{Proposition}
The involution $W = \displaystyle \frac{1}{\sqrt{2}} \begin{bmatrix} 0 & -\I_2 \\ 2 \I_2 & 0\end{bmatrix}$ acts on 
$A(\G_0(2)_2)_{even}$ as follows
\begin{align*}
 \vartheta(W \cdot \tau) = (2 \det \tau)^2 \vartheta(\tau), \qquad \phi_1(W \cdot \tau) = 1024 (2 \det \tau)^4 \phi_2(\tau ) \\
\phi_2(W \cdot \tau) = (2 \det \tau)^4 \phi_1(\tau) /1024, \qquad \chi(W \cdot \tau) = (2 \det \tau)^6 \chi(\tau).
\end{align*}
Therefore we have 
\[
A(\G_0^*(2)_2)_{even} = \CC[\vartheta, \phi, \chi, \psi], \qquad
\mathrm{Proj} \medspace A(\G_0^*(2)_2)_{even} \cong \PP(2:4:6:8)
\]
where $\phi = \phi_1 + 1024\phi_2$ and $\psi = \phi_1 \phi_2$.
\\ \\
{\bf Proof.}
By the following formula (\cite{Ig1}, p.408) 
\begin{align*}
\theta_{0000}^2(\tau/2) &= \theta_{0000}^2(\tau) + \theta_{1000}^2(\tau) + \theta_{0100}^2(\tau) + \theta_{1100}^2(\tau) \\
\theta_{0001}^2(\tau/2) &= \theta_{0000}^2(\tau) + \theta_{1000}^2(\tau) - \theta_{0100}^2(\tau) - \theta_{1100}^2(\tau) \\
\theta_{0010}^2(\tau/2) &= \theta_{0000}^2(\tau) - \theta_{1000}^2(\tau) + \theta_{0100}^2(\tau) - \theta_{1100}^2(\tau) \\
\theta_{0011}^2(\tau/2) &= \theta_{0000}^2(\tau) - \theta_{1000}^2(\tau) - \theta_{0100}^2(\tau) + \theta_{1100}^2(\tau) \\
\theta_{0100}^2(\tau/2) &= 2(\theta_{0000} \theta_{0100} + \theta_{1000} \theta_{1100})(\tau) \\
\theta_{0110}^2(\tau/2) &= 2(\theta_{0000} \theta_{0100} - \theta_{1000} \theta_{1100})(\tau) \\
\theta_{1000}^2(\tau/2) &= 2(\theta_{0000} \theta_{1000} + \theta_{0100} \theta_{1100})(\tau) \\
\theta_{1001}^2(\tau/2) &= 2(\theta_{0000} \theta_{1000} - \theta_{0100} \theta_{1100})(\tau) \\
\theta_{1100}^2(\tau/2) &= 2(\theta_{0000} \theta_{1100} + \theta_{0100} \theta_{1000})(\tau) \\
\theta_{1111}^2(\tau/2) &= 2(\theta_{0000} \theta_{1100} - \theta_{0100} \theta_{1000})(\tau)
\end{align*} 
we see that
\begin{align*}
\vartheta(W \cdot \tau) &= (\theta_{0000}^4 + \theta_{0001}^4 + \theta_{0010}^4 + \theta_{0011}^4)(-\tau^{-1}/2)/4 
= (\theta_{0000}^4 + \theta_{1000}^4 + \theta_{0100}^4 + \theta_{1100}^4)(-\tau^{-1}). 
\end{align*}
Applying the inversion formula, we obtain 
\[
 \vartheta(W \cdot \tau) = 4 (\det \tau)^2 \vartheta(\tau).
\]
By the same way, we can show that
\[
 \phi_2(-\tau^{-1}/2) = (\det \tau)^4 \phi_1(\tau) /64,
\]
and replaceing $\tau$ by $-\tau^{-1}/2$, we see that
\[
 \phi_1(-\tau^{-1} / 2) = 16384 (\det \tau)^4 \phi_2(\tau ).
\]
For the modular form $\chi(\tau)$, we have
\begin{align*}
\frac{64}{(\det \tau)^6} \chi (-\tau^{-1}/2) = (\theta_{0000}^2 \theta_{0001}^2 - \theta_{0010}^2 \theta_{0011}^2)
 (\theta_{0000}^2 \theta_{0010}^2 - \theta_{0001}^2 \theta_{0011}^2)
(\theta_{0000}^2 \theta_{0011}^2 - \theta_{0001}^2 \theta_{0010}^2)(\tau).
\end{align*}
Since the right hand side vanishes on $\HH \times \HH$, it coincides with $c \chi(\tau)$ for some constant $c$. 
Conparing Fourier coefficients, we see that $c = 1$.
\hfill $\Box$
\section{Boundary}
\subsection{}
Let us study the extension of the period map $Per : \mathcal{U} \longrightarrow \DD_{ns}^+$ to the locus 
$\{u_3 = 0\}$. Note that we have
\[
 F_C(1, \frac{1}{2} ; 1,1,1 ; -2u_1, -2u_2, 0) = F_4(1, \frac{1}{2} ; 1,1 ; -2u_1, -2u_2)
\]
where $F_4$ is Appell's hypergeometric series , and we have
\[
 F_4(1,\frac{1}{2};1,1; \frac{-x}{(1-x)(1-y)}, \frac{-y}{(1-x)(1-y)}) = 
(1-x)^{\frac{1}{2}} (1-y)^{\frac{1}{2}} {}_2F_1(\frac{1}{2}, \frac{1}{2};1;xy)
\]
(see \cite{E}). It is known that Gauss's hypergeometric series ${}_2F_1(\frac{1}{2}, \frac{1}{2};1;t)$ has an elliptic 
integral representation. In deed, the same computation in Proposition \ref{FC-prop} shows that 
$F_4(1, \frac{1}{2} ; 1,1 ; -2u_1, -2u_2)$ is a period integral of a curve
\[
 C(u) : xy(x + y + 1) + u_1 y + u_2 x = 0 
\]
of degree $(2,2)$ in $\PP^1 \times \PP^1$. The relation between this family and Appell's $F_4$ was alredy studied by Stienstra in \cite{St}. 
Here we study the relation between the invariants of $C(u)$ and the degeneration of the map
\[
 DvG : \mathcal{U} \longrightarrow \PP(1,2,3,4), \quad (u_1,u_2,u_3) \mapsto [-4s_1 + 1: s_2: 2s_3: s_1s_3]
\]
defined by invariants of cubic surfaces, where $s_i$ is the $i$-th symmetric polynomial of $u_1, u_2, u_3$. 
For this map, we have
\[
 \lim_{u_3 \rightarrow 0 } [-4s_1 + 1: s_2: 2s_3: s_1s_3] = [1-4u_1 -4u_2: u_1 u_2:0:0].
\]
On the other hand, the curve $C(u)$ is birationally equivalent to an elliptic curve
\[
 E(u) \ : \ Y^2 = f_u(X) = X^4 + X^3 + (-u_2 + \frac{u_1}{2} + \frac{1}{4}) X^2 + \frac{u_1}{4} X 
+ \frac{u_1^2}{16}
\]
by the transformation
\[
(x, y) = (2X, \frac{4Y - 4X^2 - 2 X - u_1}{4X}).
\]
\subsection{Lemma}
\label{u1-u2-inv}
The classical invariants of the quartic equation $f_u(X) = 0$ are 
\begin{align*}
 g_2(u) &= \frac{1}{192}((1- 4u_1 - 4u_2)^2 - 48 u_1 u_2), \\
g_3(u) &= -\frac{1}{13824}(1- 4u_1 - 4u_2)((1- 4u_1 - 4u_2)^2 -72 u_1 u_2), \\
\Delta_E(u) & = g_2(u)^3 - 27 g_3(u)^2 = \frac{1}{4096} u_1^2 u_2^2 ((1- 4u_1 - 4u_2)^2 -64 u_1 u_2). 
\end{align*}
Therefore $[1- 4u_1 - 4u_2:u_1 u_2] \in \PP(1,2)$ corresponds to a singular $E(u)$ iff 
\[
 [1- 4u_1 - 4u_2:u_1 u_2] = [1:0] \ \text{or} \ [8:1].
\]
Moreover we have $\Delta_{sing}(u_1, u_2, 0) = (4096 \Delta_E(u_1, u_2) / u_1^2 u_2^2)^2$.
\\ \\
{\bf Proof.} 
This is obtained from the definition
\begin{align*}
g_2 &= ae - 4bd + 3c^2, \quad 
g_3 = \det \begin{bmatrix} a & b & c \\ b & c & d \\ c & d & e \end{bmatrix}
\end{align*}
for $aX^4 + 4 bX^3 + 6cX^2 + 4dX + e = 0$.
\hfill $\Box$
\subsection{}
Now we can define a degenerated period map
\[
 Per_{12} : \mathcal{U}_{12} = \{(u_1, u_2) \in \CC^2 \ | \ \Delta_E(u_1,u_2) \ne 0 \} 
\longrightarrow \HH,
\]
and construct the inverse map 
\[
 \HH \longrightarrow [1- 4u_1 - 4u_2:u_1 u_2] \in \PP(1,2)
\]
by the Siegel $\Phi$-operator $\displaystyle \Phi(f)(\tau_1) = \lim_{t \rightarrow \infty} f(\begin{bmatrix} \tau_1 & 0 \\ 0 & it\end{bmatrix})$. 
Let us define modular forms
\begin{align*}
h_1 &= \Phi(8 \vartheta) = 4(\theta_{00}^4 + \theta_{01}^4), \\
h_2 &= \Phi(\vartheta^2 - \phi) = 
\frac{1}{4}(\theta_{00}^4 + \theta_{01}^4)^2 - \theta_{00}^4 \theta_{01}^4 =
\frac{1}{4}(\theta_{00}^4 - \theta_{01}^4)^2
\end{align*} 
of weight $2$ and $4$ with respect to $\G_0(2)_1$.  
\subsection{Lemma}
Modular forms $h_1$ and $h_2$ satisfy same relations for $1-4u_1-4u_2$ and $u_1 u_2$ in Lemma \ref{u1-u2-inv}. In deed, we have
\begin{align*}
h_1(\tau)^2 - 48 h_2(\tau) &= 64 E_4(2\tau) = 64(1 + 240 \sum_{n=1}^{\infty} \sigma_3(n)q^{2n}), \\
h_1(\tau)(h_1(\tau)^2 - 72 h_2(\tau)) &= -512E_6(2\tau) = -512(1 - 504 \sum_{n=1}^{\infty} \sigma_5(n)q^{2n}), \\
h_2(\tau)^2 (h_1(\tau)^2 - 64 h_2(\tau)) &= 2^{18} \eta(2 \tau) = 2^{18} q^2 \prod_{n=1}^{\infty}(1-q^{2n})^{24}, \\
h_2(\tau) / (h_1(\tau)^2 - 64 h_2(\tau)) &= \eta(2\tau) / \eta(\tau)
\end{align*}
where $q = \exp(2 \pi i \tau)$, and
\begin{align*}
\lim_{t \rightarrow \infty} [h_1(it):h_2(it)] = [1:0], \quad 
\lim_{t \rightarrow \infty} [h_1(-1/2it):h_2(-1/2it)] = [8:1] \in \PP(2,4).
\end{align*}
Since $\eta(2\tau)/ \eta(\tau)$ is the Hauptmodul for $\G_0(2)_1$, we see that the map
\[
 \HH / \G_0(2)_1 \cup \{0, \infty\} \rightarrow \PP(1,2), \quad \tau \mapsto [h_1(\tau): h_2(\tau)] = [1- 4u_1 - 4u_2:u_1 u_2] 
\]
is an isomorphism.
\\ \\
{\bf Proof.}
By the formula
\begin{align*}
\theta_{00}^2(2 \tau) = \frac{1}{2}(\theta_{00}^2(\tau) + \theta_{01}^2(\tau)), \quad
\theta_{01}^2(2 \tau) = \theta_{00}(\tau) \theta_{01}(\tau), \quad
\theta_{10}^2(2 \tau) = \frac{1}{2}(\theta_{00}^2(\tau) - \theta_{01}^2(\tau)),
\end{align*}
we have
\begin{align*}
E_4(2\tau) &=  [\theta_{00}^8 - \theta_{00}^4 \theta_{01}^4 + \theta_{01}^8](2\tau) \\
&= \frac{1}{16}[(\theta_{00}^2 + \theta_{01}^2)^4 - 4(\theta_{00}^2 + \theta_{01}^2)^2 \theta_{00}^2 \theta_{01}^2 
+ 16\theta_{00}^4 \theta_{01}^4](\tau) \\
& = \frac{1}{64}[h_1^2 - 48 h_2](\tau)
\end{align*}
and
\begin{align*}
E_6(2\tau) &= -\frac{1}{2}[(\theta_{00}^4 + \theta_{01}^4) (2\theta_{00}^4 - \theta_{01}^4)(\theta_{00}^4 - 2\theta_{01}^4)](2\tau) \\
&= -\frac{1}{64}[(\theta_{00}^4 +6 \theta_{00}^2 \theta_{01}^2 + \theta_{01}^4) 
(\theta_{00}^4 + \theta_{01}^4)(\theta_{00}^4 -6 \theta_{00}^2 \theta_{01}^2 + \theta_{01}^4)](\tau) \\
&= -\frac{1}{512}[h_1 (h_1^2 - 72 h_2)](\tau).
\end{align*}
Other assertions are shown by similar calculation.
\subsection{Theorem}
Let us define an embedding $\Theta : \Si_2 / \G_0^*(2)_2 \rightarrow \PP(1,2,3,4)$ by
\[
\tau \mapsto [8 \vartheta:\vartheta^2 - \phi :1024 \chi:1024(\psi - \vartheta \chi)]
\]
Then we have the commutative diagram
\[
\begin{CD}
\mathcal{U} @>Per>> \DD_{ns}^+ \\
@V DvG VV @VVV \\
\PP(1,2,3,4) @<\Theta<< \Si / \G_0^*(2)_2 \cong \DD_{ns}^+ / \Ons
\end{CD}
\]
and $\Theta$ induces an isomorphism $\Si_2 / \G_0^*(2)_2 \cup \HH / \G_0(2)_1 \cup \{0, \infty\} \cong \PP(1,2,3,4)$.
\\ \\ 
{\bf Proof.}  
In deed, the map $\Theta$ is the unique map 
\[
 \Si / \G_0^*(2)_2 \rightarrow \PP(1,2,3,4), \quad \tau \mapsto 
[F_2(\tau):F_4(\tau):F_6(\tau):F_8(\tau)] \qquad (F_k \in M_k(\G_0^*(2)_2))
\]
such that 
\\
(i) $F_6$ vanishes on $\HH \times \HH$, 
\\
(ii) 
\[
 \lim_{t \rightarrow \infty} [F_2:F_4:F_6:F_8](\begin{bmatrix} \tau & 0 \\ 0 & i\infty \end{bmatrix}) = [h_1(\tau):h_2(\tau):0:0],
\]
(iii) 
\begin{align*}
[F_2:F_4:F_6:F_8](\begin{bmatrix} 2\tau & \tau \\ \tau & 2\tau \end{bmatrix}) &= 
[-4s_1 + 1: s_2: 2s_3: s_1s_3]|_{u_1 = u_2 = u_3=u} \\
&=[-12u +1:3u^2:2u^3:3u^4],
\end{align*}
that is, $F_4^2 = 3 F_8$ and $9F_6^2 = 4F_4 F_8$.
\\
For (i), we see that $F_6 = c \chi$ by Proposition \ref{prop-table} and Lemma \ref{chi-vanish}. For (ii), note that
\begin{align*}
 &\lim_{t \rightarrow \infty} (c_1 \vartheta^4 + c_2 \vartheta^2 \phi + c_3 \phi^2 + c_4 \vartheta \chi + c_5 \psi)
(\begin{bmatrix} \tau_1 & 0 \\ 0 & it \end{bmatrix}) = 0\\
\Leftrightarrow \ &c_1 (\theta_{00}^4 + \theta_{01}^4)^4/16 + c_2 (\theta_{00}^4 + \theta_{01}^4)^2 \theta_{00}^4 \theta_{01}^4/4
+ c_3 (\theta_{00}^4 \theta_{01}^4)^2 =0 \\
\Leftrightarrow \ &c_1 = c_2 = c_3 = 0
\end{align*}
Therefore we have
\[
 [F_2:F_4:F_6:F_8] = [8 \vartheta:\vartheta^2 - \phi :c \chi:c_4 \vartheta \chi + c_5 \psi].
\]
Now the condition (iii) implies
\[
 (\vartheta^2 - \phi)^2(\tau) = 3(c_4 \vartheta \chi + c_5 \psi)(\tau), \quad 
9(c \chi)^2(\tau) = [4(\vartheta^2 - \phi)(c_4 \vartheta \chi + c_5 \psi)](\tau)
\]
for $\tau = \begin{bmatrix} 2\tau_1 & \tau_1 \\ \tau_1 & 2\tau_1 \end{bmatrix}$. Conparing Fourier coefficients 
\begin{align*}
\vartheta(\tau) &= 1 + 72 q^8 + 192 q^{12} + 504 q^{16} + 576 q^{20} + 2280 q^{24} + \cdots, \\
\phi(\tau) &= q^8 - 4q^{12} - 2 q^{16} + 20 q^{20} + 5 q^{24} + \cdots, \\
\chi(\tau) &= q^{12} - 6 q^{16} + 3 q^{20} + 40 q^{24} + \cdots, \\
\psi(\tau) &= q^{12} + 6 q^{16} -21 q^{20} -56 q^{24} + \cdots,
\end{align*}
we cobtain $c = 1024, c_4 = -1024$ and $c_5 = 1024$.


\begin{thebibliography}{20}
\bibitem[B]{B}
 V. V. Batyrev, {\it Dual polyhedra and mirror symmetry for Calabi-Yau 
        hypersurfaces in toric varieties}, J. Alg. Geom, 3 (1994), 493-535.

\bibitem[C1]{C1}
 D. A. Cox, {\it Minicourse on Toric Varieies}, Lecture note (2001).

\bibitem[C2]{C2}
 D. A. Cox, {\it Primes of the form $x^2+ny^2$}, John Wiley $\&$ Sons Inc., New York (1989).

\bibitem[D1]{D1}
 I. V. Dolgachev, {\it Mirror Symmetry for Lattice Polarized K3 Surfaces}, J. Math. Sci., vol. 81 (1996), No. 3, 
         2599-2630.

\bibitem[D2]{D2}
 I. V. Dolgachev, {\it Lectures on Modular Forms}, Lecture note (2005).

\bibitem[DK]{DK}
 I. Dolgachev and J. Keum, {\it Birational automorphisms of quartic Hessian surfaces},
        Trans. Amer. Math. Soc. 354 (2002), 3031-3057.

\bibitem[DvG]{DvG}
 E. Dardanelli and B. van Geemen, {\it Hessians and the moduli space of cubic surfaces},
        Contemp. Math. 422, AMS (2007) 17-36.

\bibitem[E]{E}
 A. Erdelyi ed., Higher Transcendental Functions, Bateman Manuscript Project (McGraw-Hill, New York, 1953). 

\bibitem[GN]{GN}
 V. A. Gritsenko and V. V. Nikulin, {\it K3 surfaces, Lorentzian Kac-Moody Algebras and mirror symmetry}, 
        Math. Res. Let., 3 (1996), 211-229. 

\bibitem[H]{H}
 B. Hunt, {\it The Geometry of some special Arithmetic Quotients}, 
 Springer LNM 1637 (1996).
 
\bibitem[Ib]{Ib}
 T. Ibukiyama, {\it On Siegel modular varieties of level 3},  Int. J. Math. 2 (1991), No.1, 17-35.

\bibitem[Ig1]{Ig1}
 J. Igusa, {\it On Siegel modular forms of genus two. II}, Amer. J. Math. 86 (1964), 392-412. 

\bibitem[Ig2]{Ig2}
 J. Igusa, {\it Theta functions}, Springer-Verlag (1972).

\bibitem[Kb]{Kb}
 N. Koblitz, {\it Introduction to Elliptic Curves and Modular Forms}, GTM. No. 97, Springer-Verlag (1984).

\bibitem[Kl]{Kl}
 H. Klingen, {\it Introductory lectures on Siegel modular forms}, Cambridge Univ. Press (1990).

\bibitem[Ko1]{Ko1}
 K. Koike, {\it K3 surfaces induced from polytopes} (in Japanese), Master thesis, Chiba Univ. (1998).

\bibitem[Ko2]{Ko2}
 K. Koike, {\it Moduli space of Hessian K3 surfaces and arithmetic quotients}, arXiv:math 1002.2854v1.

\bibitem[M]{M}
 D. Mumford, {\it Tata Lectures on Theta I}, Birkh\"auser (1983).

\bibitem[N]{N}
 V. V. Nikulin, {\it Integral symmetric bilinear forms and some of their
	applications}, Math. USSR Izv. 14 (1980), 103-166. 

\bibitem[P]{P}
C. Peters, {\it Monodromy and Picard-Fuchs equations for families of $K3$-surfaces and elliptic curves}, 
        Ann. sci. \'Ecole Norm. Sup., Ser. (4), 19, no. 4 (1986), 583-607. 

\bibitem[PS]{PS}
 C. Peters and J. Stienstra, {\it A pencil of K3-surfaces related to Ap\'ery's recurrence for $\zeta(3)$ and 
        Fermi surfaces for potential zero}, Arithmetic of complex manifolds, LNM 1399 (1989), 110-127. 

\bibitem[Sa]{Sa} 
 G. Salmon, {\it A Treatise on the Analytic Geometry of Three Dimensions}, 
       Chelsea Publishing Company, New York (1965).

\bibitem[St]{St}
 J. Stienstra, {\it A variation of mixed Hodge structure for a special case of Appell's $F_4$}, 
       Special Differential Equations, Proceedings of the Taniguchi Workshop 1991, p.109-116. 

\bibitem[V]{V} 
 H. A. Verrill, {\it The $L$-series of certain rigid Calabi- Yau threefolds}, 
       J. Number Theory, 81 (2000), 310-334.

\bibitem[Y]{Y}
 M. Yoshida, {\it Fuchsian differential equations, Aspects of Mathematics}, 
        E11. Friedr. Vieweg \& Sohn, Braunschweig(1987). 
\end{thebibliography}
\end{document}